\documentclass[11pt,a4paper]{article}
\usepackage{amsfonts}
\usepackage{epic}
\usepackage{graphicx}
\usepackage{flafter}
\usepackage{geometry}
\usepackage[centertags]{amsmath}
\usepackage{amsfonts}
\usepackage{amssymb}
\usepackage{dsfont}
\usepackage{amsthm}
\usepackage{bbm}
\usepackage{mathrsfs}
\usepackage{epsfig}
\usepackage{color}
\usepackage{amsmath}
\numberwithin{equation}{section}
\numberwithin{figure}{section}

\date{}
\DeclareMathAlphabet{\mathds}{OT1}{cmss}{m}{sl}
\begin{document}

\title{\bf{Strongly Independent Matrices and Applications on the Rigidity of $A$-Invariant Measures on $n$-Torus}}

\date{\today}
\author{Huichi Huang, Enhui Shi,  Hui Xu\footnote{Corresponding author}}
\maketitle
$$\text{\textbf{Abstract}}$$
\indent   We introduce the notion of strongly independent matrices and show the existence of strongly independent matrices in $GL(n,\mathbb{Z})$ over $\mathbb{Z}^n\setminus\{0\}$ when $2n+1$ is a prime number. As an application of strong independence,  we give a measure rigidity result for endomorphisms on $n$-torus $\mathbb{T}^n$.\\[1cm]
{\bf Keywords}: Strong independence, ergodicity, mixing, Fourier coefficient, measure rigidity.
\newtheorem{theorem}{Theorem}[section]
\newtheorem{prop}{Proposition}
\newtheorem{definition}{Definition}[section]
\section{Introduction and main results}

For an integer $n$, consider the automorphism  $T_n$  on $\mathbb{T}=\mathbb{R}/\mathbb{Z}$ (called $\times n$ map) given by
\begin{displaymath}
T_n(x)=nx\mod 1,~x\in\mathbb{T}.
\end{displaymath}
\indent In \cite{F} H. Furstenberg showed that any closed subset of $\mathbb{T}$ invariant under the action of a non-lacunary semigroup of integers must be either finite or the whole $\mathbb{T}$. Here non-lacunary means that not to be all powers of one integer.
Then he conjectured a similar statement holds for invariant measures.
\newtheorem*{con}{Furstenberg's Conjecture}
\begin{con}
An ergodic invariant  measure on $\mathbb{T}$  under the action of a non-lacunary semigroup of integers is either Lebesgue measure or finitely supported.
\end{con}
In general, consider  two commuting maps $\times p$ and $\times q$, where $p$ and $q$ are two positive integers satisfying $\frac{\log p}{\log q}\notin \mathbb{Q}$.
In this direction, the first substantial result was obtained by R. Lyons.
\begin{theorem}[\cite{L}]
Suppose $p,q$ are two positive integers with $\frac{\log p}{\log q}\notin \mathbb{Q}$. Then any $\times p,\times q$-invariant ergodic  measure $\mu$ is either Lebesgue measure or finitely supported if $(\mathbb{T}, \mathscr{B},\mu, T_q)$ has no zero entropy factor.
\end{theorem}
This result was improved by D. Rudolph  with the assumption that $p$ and $q$ are relatively prime.
\begin{theorem}[\cite{R}]
Let $p$ and $q$ are relatively prime positive integers. Then any $\times p,\times q$-invariant ergodic measure is either Lebesgue measure or finitely supported
if one of $T_p$ and $T_q$ has positive entropy.
\end{theorem}

One may consult
\cite{KK, KS1, KS2} for the extensions of above results to automorphisms on $n$-torus with $n\geq 2$.
Recently, H. Huang has gotten the following rigidity theorem assuming that the measure is invariant under enough $\times q$ maps.
\begin{theorem}[\cite{H}]
The Lebesgue measure is the only non-atomic $\times p$-invariant measure on $\mathbb{T}$ satisfying one of the following:
\begin{itemize}
  \item [(1)] it is ergodic and there exist a nonzero integer $l$ and a F{\o}lner sequence $\Sigma=\{F_n\}_{n=1}^\infty$ in $\mathbb{N}$
  such that $\mu$  is $\times(p^j+l)$-invariant for all $j$ in some $E\subseteq \mathbb{N}$ with $D_{\Sigma}(E)=1$;
  \item [(2)]    it is weakly mixing and there exist a nonzero integer $l$ and a F{\o}lner sequence $\Sigma=\{F_n\}_{n=1}^\infty$ in $\mathbb{N}$
  such that $\mu$  is $\times(p^j+l)$-invariant for all $j$ in some $E\subseteq \mathbb{N}$ with $\overline{D}_{\Sigma}(E)>0$;
  \item [(3)]  it is strongly mixing and there exist a nonzero integer $l$ and an infinite set $E\subseteq \mathbb{N}$
  such that $\mu$  is $\times(p^j+l)$-invariant for all $j$ in  $E$.
\end{itemize}
Moreover, a $\times p$- invariant measure satisfying (2) or (3) is either a Dirac measure or the Lebesgue measure.
\end{theorem}

In this paper, we extend the above measure rigidity theorem to endomorphisms on $n$-torus $\mathbb T^n$. Here we use the notion of strongly independent matrices (see Definition \ref{independent}).

\begin{theorem}\label{rigidity on torus}
The Lebesgue measure is the only non-atomic $\times A$-invariant measure on $\mathbb{T}^n$ satisfying one of the following:
\begin{itemize}
  \item [(1)] it is ergodic and there exist a sequence of $n$ matrices $B_1, B_2, \cdots, B_n\in GL(n,\mathbb{Z})$ which are strongly independent over $\mathbb{Z}^n\setminus\{0\}$ and a F{\o}lner sequence
   $\Sigma=\{F_n\}_{n=1}^\infty$ in $\mathbb{N}$ such that $\mu$  is $\times(A^j+B_i)$-invariant for all $j$ in some $E\subseteq \mathbb{N}$ with $D_{\Sigma}(E)=1$ and all $i=1,2,\cdots,n$;
  \item [(2)]    it is weakly mixing and there exist a sequence of $n$ matrices $B_1, B_2, \cdots, B_n\in GL(n,\mathbb{Z})$ which are strongly independent over $\mathbb{Z}^n\setminus\{0\}$ and a F{\o}lner sequence $\Sigma=\{F_n\}_{n=1}^\infty$ in $\mathbb{N}$ such that $\mu$  is $\times(A^j+B_i)$-invariant for all $j$ in some $E\subseteq \mathbb{N}$ with $\overline{D}_{\Sigma}(E)>0$ and all $i=1,2,\cdots,n$;
  \item [(3)]  it is strongly mixing and there exist a sequence of $n$ matrices $B_1, B_2, \cdots, B_n\in GL(n,\mathbb{Z})$ which are strongly independent over $\mathbb{Z}^n\setminus\{0\}$  and an infinite set $E\subseteq \mathbb{N}$
  such that $\mu$  is $\times(A^j+B_i)$-invariant for all $j$ in  $E$ and all $i=1,2,\cdots,n$.
\end{itemize}
Moreover, a $\times A$-invariant measure satisfying (2) or (3) is either a Dirac measure or the Lebesgue measure.
\end{theorem}
In the following theorem, we give some concrete examples of strongly independent matrices. These examples show that one can construct a countably generated abelian semigroup action on $\mathbb{T}^n$ such that the Lebesgue measure is the only non-atomic invariant measure under the action.
For the definition of $M_n(2)$, please see Section 2.

\begin{theorem}\label{strongly independent matrix}
If $2n+1$ is a prime number then the matrices $M_n(2),M_n^2(2),\cdots,M_n^n(2)$ are strongly independent over $\mathbb{Z}^n\setminus\{0\}$.
\end{theorem}
The paper is organized as follows. In section 2, we give some definitions and notations which will be used throughout the paper. In section 3, we give some lemmas about cyclotomic polynomial and some basic facts about Fourier coefficients. We obtain a criterion to determine the strong independence of the powers of a matrix and construct an infinite family of strongly independent matrices in Section 4.
We  characterize the ergodicity, weakly mixing, and strongly mixing by Fourier coefficients in section 5. In the last section, we prove a measure rigidity theorem for endomorphisms on $\mathbb T^n$.

\section{Preliminaries}
Let $\mathbb{N}$ stand for the set of nonnegative integers and $|E|$ denote the cardinality of a set $E$. Denote by $M(n,\mathbb{Z})$ the set of square matrices of order $n$ with entries in $\mathbb{Z}$. Denote by $GL(n,\mathbb{Z})$  the set of invertible matrices in $M(n,\mathbb{Z})$ whose inverse is also in $M(n,\mathbb{Z})$ .\\

Let $M_n(a)$ and $N_n(a)$ denote the following $n\times n$ real tridiagonal matrices with $n$ being a positive integer and $a$ being a real number.
\begin{displaymath}
M_n(a) =\begin{bmatrix}
a&1&0&\cdots&0&0\\
1&a&1&\cdots&0&0\\
0&1&a&\cdots&0&0\\
\vdots&\vdots&\vdots&\cdots&\vdots&\vdots\\
0&0&0&\cdots&a&1\\
0&0&0&\cdots&1&a-1
\end{bmatrix},~~~~~~~
N_n(a) =\begin{bmatrix}
a&1&0&\cdots&0&0\\
1&a&1&\cdots&0&0\\
0&1&a&\cdots&0&0\\
\vdots&\vdots&\vdots&\cdots&\vdots&\vdots\\
0&0&0&\cdots&a&1\\
0&0&0&\cdots&1&a
\end{bmatrix}.
\end{displaymath}
In particular,
\begin{displaymath}
M_1(a) =\begin{bmatrix}
a-1
\end{bmatrix},~~~~
M_2(a) =\begin{bmatrix}
a&1\\1&a-1
\end{bmatrix},~~~~
N_1(a) =\begin{bmatrix}
a
\end{bmatrix},~~~~
N_2(a) =\begin{bmatrix}
a&1\\1&a
\end{bmatrix}.
\end{displaymath}
\begin{definition}\label{independent}
We call  $n$ matrices $B_1, B_2, \cdots, B_n\in GL(n,\mathbb{Z})$ \textbf{strongly independent} over $\mathbb{Z}^n\setminus\{0\}$ if for any nonzero row vector $\vec{k}=(k_1,k_2,\cdots,k_n)\in\mathbb{Z}^n$, the vectors $\vec{k}B_1,\vec{k}B_2,\cdots,\vec{k}B_n$ are independent over $\mathbb{R}$.
\end{definition}

\begin{definition}
A \textbf{ F{\o}lner sequence} in $\mathbb{N}$ is a sequence $\Sigma=\{F_n\}_{n=1}^\infty$ of finite subsets in $\mathbb{N}$ satisfying
\begin{displaymath} \lim_{n\rightarrow \infty}\frac{|(F_n+m)\Delta F_n|}{|F_n|}=0 \end{displaymath}
for every $m$ in $\mathbb{N}$. Here $\Delta$ stands for the symmetric difference.
\end{definition}
\begin{definition}
Let $\Sigma=\{F_n\}_{n=1}^\infty$ be a sequence of finite subsets of $\mathbb{N}$. The \textbf{density} $D_{\Sigma}(E)$ of a subset $E$ of
$\mathbb{N}$ with respect to $\Sigma$ is given by \begin{displaymath} D_{\Sigma}(E):=\lim_{n\rightarrow \infty}\frac{|E\cap F_n|}{|F_n|}. \end{displaymath}
The upper density $ \overline{D}_{\Sigma}(E)$ and lower density $ \underline{D}_{\Sigma}(E)$ are given by
\begin{displaymath} \overline{D}_{\Sigma}(E):=\limsup_{n\rightarrow \infty}\frac{|E\cap F_n|}{|F_n|} , ~~~~~  \underline{D}_{\Sigma}(E):=\liminf_{n\rightarrow \infty}\frac{|E\cap F_n|}{|F_n|}.\end{displaymath}
\end{definition}

\begin{definition}
For $\vec{k}=(k_1,k_2,\cdots,k_n)\in\mathbb{Z}^n$, the \textbf{Fourier coefficient} $\hat{\mu}(\vec{k})$ of a measure $\mu$ on $\mathbb{T}^n$ is defined by
\begin{displaymath} \hat{\mu}(\vec{k})=\int_{\mathbb{T}^n} z_1^{k_1}z_2^{k_2}\cdots z_n^{k_n}d\mu(z_1,z_2,\cdots,z_n)\end{displaymath} when taking $\mathbb{T}^n=\{(z_1,z_2,\cdots,z_n)\in\mathbb{C}^n:~|z_1|=|z_2|=\cdots=|z_n|=1\}$.
\end{definition}
Within this paper, a measure on a compact metrizable $X$ always means a Borel probability measure. A measure $\mu$ is called \textbf{non-atomic} if $\mu(\{x\})=0$ for every $x\in X$.\\
\indent A measure $\mu$ on $X$ is called $T$-\textbf{invariant} for a continuous map $T: X\rightarrow X$ if $\mu(E)=\mu(T^{-1}E)$ for every Borel subset $E$ in $X$. A $T$-invariant measure $\mu$ is called \textbf{ergodic} if every Borel subset $E$ with $T^{-1}E=E$ implies $\mu(E)=0$ or $1$, it is called \textbf{weakly mixing} if $\mu\times \mu$ is an ergodic $T\times T$-invariant measure on $X\times X$, and it is called \textbf{strongly mixing} if $\lim_{j\rightarrow \infty}\mu(T^{-j}E\cap F)=\mu(E)\mu(F)$ for all Borel subsets $E,F$ in $X$.\\
\indent Within this paper, we only consider that $X=\mathbb{T}^n=\mathbb{R}^n/\mathbb{Z}^n$ and $T=T_A$ is the $\times A$ map on  $\mathbb{T}^n$ defined by
$ T_A(x)=A x\mod \mathbb{Z}^n$ for all column vectors $x\in \mathbb{R}^n/\mathbb{Z}^n$ and $A\in M(n,\mathbb{Z})$. \\
\indent For $\vec{k}=(k_1,k_2,\cdots,k_n)\in\mathbb{Z}^n$ and $z=(z_1,z_2,\cdots,z_n)\in \mathbb{T}^n$,
\begin{displaymath} z^{\vec{k}}:=z_1^{k_1}z_2^{k_2}\cdots z_n^{k_n}.\end{displaymath}

\section{Some Lemmas}
\newtheorem{lemma}{Lemma}[section]
Since $n^{\text{nt}}$ cyclotomic polynomial is irreducible in $\mathbb{Q}[X]$, we have the following result which also can be found in \cite{W}.
\begin{lemma}\label{cycl1}
Let $n$ be a positive integer and $\zeta_n=e^{\frac{2\pi i}{n}}$. Then
\begin{displaymath}
\left[\mathbb{Q}(\zeta_n):\mathbb{Q} \right]=\varphi(n),
\end{displaymath}
where $\varphi$ is Euler function.
\end{lemma}

\begin{lemma}\label{cycl2}
Let $n$ be a positive integer and $\zeta_{2n+1}=e^{\frac{2\pi i}{2n+1}}$. Then
\begin{itemize}
  \item [(1)] \begin{math}\left[\mathbb{Q}(\zeta_{2n+1}):\mathbb{Q}(\cos\frac{2\pi}{2n+1}) \right]= 2. \end{math}
  \item [(2)] \begin{math}\left[\mathbb{Q}(\zeta_{2n+1}):\mathbb{Q}(\cos\frac{2\pi}{2n+1},\cos\frac{4\pi}{2n+1},\cdots,\cos\frac{2n\pi}{2n+1}) \right]= 2. \end{math}
\end{itemize}
\end{lemma}

\noindent \textit{Proof}. (1) Since $\zeta_{2n+1}^2-2\left(\cos\frac{2\pi}{2n+1}\right)\zeta_{2n+1}+1=0$ and $x^2-2\left(\cos\frac{2\pi}{2n+1}\right)x+1=0$ is irreducible over $\mathbb{Q}\left(\cos\frac{2\pi}{2n+1}\right)$, we have
 \begin{displaymath}\left[\mathbb{Q}(\zeta_{2n+1}):\mathbb{Q}\left(\cos\frac{2\pi}{2n+1}\right) \right]= 2. \end{displaymath}
(2) Noting that $\cos(k+1)\theta=2\cos k\theta\cos\theta-\cos(k-1)\theta$, we can derive inductively that all $\cos\frac{2\pi}{2n+1},\cos\frac{4\pi}{2n+1},\cdots,\cos\frac{2n\pi}{2n+1}$ are in  $\mathbb{Q}(\cos\frac{2\pi}{2n+1})$, then
\begin{displaymath}
\left[\mathbb{Q}\left(\cos\frac{2\pi}{2n+1},\cos\frac{4\pi}{2n+1},\cdots,\cos\frac{2n\pi}{2n+1}\right): \mathbb{Q}\left(\cos\frac{2\pi}{2n+1}\right)\right]= 1.
\end{displaymath}
The result holds by
\begin{eqnarray*}
&&\left[\mathbb{Q}(\zeta_{2n+1}):\mathbb{Q}\left(\cos\frac{2\pi}{2n+1}\right)\right]=\\&&\left[\mathbb{Q}(\zeta_{2n+1}):\mathbb{Q}
\left(\cos\frac{2\pi}{2n+1},\cdots,\cos\frac{2n\pi}{2n+1}\right) \right]\cdot
\left[\mathbb{Q}\left(\cos\frac{2\pi}{2n+1},\cdots,\cos\frac{2n\pi}{2n+1}\right): \mathbb{Q}\left(\cos\frac{2\pi}{2n+1}\right)\right].
\end{eqnarray*}
\rightline{$\Box$}

\begin{lemma}\label{supp}
For a nonzero $\vec{k}=(k_1,k_2,\cdots,k_n)\in\mathbb{Z}^n$, if $\hat{\mu}(\vec{k})=1$ then the support of $\mu$
\begin{displaymath} \emph{Supp}(\mu)\subseteq \left\{(z_1,z_2,\cdots,z_n)\in\mathbb{T}^n:~z_1^{k_1}z_2^{k_2}\cdots z_n^{k_n}=1\right\}\end{displaymath}
\end{lemma}

\begin{lemma}\label{atomic}
Let $B_1, B_2, \cdots, B_n\in GL(n,\mathbb{Z})$ be strongly independent over $\mathbb{Z}\setminus\{0\}$.  If there is  some nonzero $\vec{k}\in\mathbb{Z}^n$ such that $\hat{\mu}({\vec{k}B_i})=1$ for every $i\in\{1,2,\cdots,n\}$, then $\mu$ is an atomic measure on $\mathbb{T}^n$.
\end{lemma}
\noindent\textit{Proof}. Let $\vec{k}B_i=(l_{1i},l_{2i},\cdots, l_{ni})$ for each $i\in\{1,2,\cdots,n\}$. By Lemma \ref{supp},
\begin{displaymath}
\text{Supp}(\mu)\subseteq \bigcap_{i=1}^n\left\{(z_1,z_2,\cdots,z_n)\in\mathbb{T}^n:~z_1^{l_{1i}}z_2^{l_{2i}}\cdots z_n^{l_{ni}}=1\right\}.
\end{displaymath}
Writing it in addition notation that is
\begin{eqnarray*}
&&\bigcap_{i=1}^n\left\{(x_1,x_2,\cdots,x_n)\in [0,1)^n:~l_{1i}x_1+l_{2i}x_2+\cdots +l_{ni}x_n\in \mathbb{Z}\right\}\\
&&\subseteq \left\{\vec{k}L^{-1}:~~\vec{k}\in [-M,M]^n\cap\mathbb{Z}^n\right\},
\end{eqnarray*}
where $M=\max\{|l_{i1}|+|l_{i2}|+\cdots+|l_{in}|:~~i=1,2,\cdots,n\}$ and $L=(l_{ij})_{n\times n}$ are constant . Therefore the support of $\mu$ is a finite set since that $L$ is invertible over $\mathbb{R}$ according to the independence of $\vec{k}B_1,\vec{k}B_2,\cdots,\vec{k}B_n$.\\

The following lemmas are stated in \cite{H}.
\begin{lemma}\label{ergodicmean}
For a topological dynamical system $(X,T)$, if $\nu$ is an ergodic $T$-invariant measure on $X$, then for every F{\o}lner sequence $\{F_n\}_{n=1}^\infty$ in $\mathbb{N}$ ,one has
\begin{displaymath} \lim_{n\rightarrow \infty} \frac{1}{|F_n|}\sum_{j\in F_n}F(T^j x)=\int_X f d\nu\end{displaymath}
for every $f\in L^(X,\nu)$(note that the identity holds with respect to $L^2$-norm). Consequently
\begin{equation}
\lim_{n\rightarrow \infty} \frac{1}{|F_n|}\sum_{j\in F_n}\int_X f(T^j x)g(x) d\nu(x)=\int_X f d\nu \int_X g d \nu\label{eq ergodic mean}
\end{equation}
for every $f,g$ in $L^2(X,\nu)$.
\end{lemma}

\begin{lemma}\label{dirac}
Let $T:X\rightarrow X$ be a continuous map on a compact metrizable space $X$. If a weakly mixing $T$-invariant atomic measure, then $\mu$ is a Dirac measure on $X$.
\end{lemma}

\section{Existence of Strongly Independent Matrices}
\begin{theorem}
If a matrix $B\in GL(n,\mathbb{Z})$ satisfies that all its eigenvalues are different reals and the vector space spanned by   any $n-1$ row eigenvectors of which doesn't contain any nonzero elements of $\mathbb{Z}^n$, then $B,B^2,\cdots,B^n$ are strongly independent over $\mathbb{Z}^n\setminus\{0\}$.
\end{theorem}
\noindent \textit{Proof}. We assume that the eigenvalues of $B$ are $\lambda_1<\lambda_2<\cdots<\lambda_n$. Let $\Lambda$ denote diag$(\lambda_1,\lambda_2,\cdots,\lambda_n)$. Then there exists some invertible real matrix $P$ such that $B=P^{-1}\Lambda P$.\\
Write $P=\begin{bmatrix}\beta_1\\ \beta_2\\ \vdots\\ \beta_n\end{bmatrix}$. Then each $\beta_i$ is an row eigenvector of $B$ corresponding to $\lambda_i$ according to $PB=\Lambda P$.\\
We claim that, for any nonzero vector $\vec{k}=(k_1,k_2,\cdots,k_n)\in\mathbb{Z}^n$, each entry of $\vec{k}P^{-1}=(x_1,x_2,\cdots,x_n)$ is nonzero. If not, we assume $x_n=0$, then $\vec{k}=(x_1,x_2,\cdots,x_{n-1},0)P=x_1\beta_1+\cdots+x_{n-1}\beta_{n-1}\in $ span$\{\beta_1,\cdots,\beta_{n-1}\}$, which contradicts to the hypothesis. Therefore
\begin{equation*}
\begin{bmatrix}
\vec{k}B\\  \vec{k}B^2\\  \vdots\\  \vec{k}B^n
\end{bmatrix}
=
\begin{bmatrix}
\vec{k}P^{-1}\Lambda P\\  \vec{k}P^{-1}\Lambda^2 P\\  \vdots\\  \vec{k}P^{-1}\Lambda^n P
\end{bmatrix}
=
\begin{bmatrix}
(x_1,x_2,\cdots,x_n)\Lambda P\\  (x_1,x_2,\cdots,x_n)\Lambda^2 P\\  \vdots\\  (x_1,x_2,\cdots,x_n)\Lambda^n P
\end{bmatrix}
=
\begin{bmatrix}
\lambda_1& \lambda_1^2&\cdots&\lambda_1^n\\
\lambda_2& \lambda_2^2&\cdots&\lambda_2^n\\
\vdots&\vdots&&\vdots\\
\lambda_n& \lambda_n^2&\cdots&\lambda_n^n\\
\end{bmatrix}
\begin{bmatrix}
x_1& & & \\
&x_2& & \\
&&\ddots&\\
&&&x_n
\end{bmatrix}
P
\end{equation*}
which is invertible over $\mathbb{R}$ since $\lambda_1,\lambda_2,\cdots,\lambda_n$ are all nonzero and different and $x_1,x_2,\cdots,x_n$
are all nonzero. Consequently, for any nonzero $\vec{k}\in\mathbb{Z}^n$, $\vec{k}B,  \vec{k}B^2,  \cdots,  \vec{k}B^n$ are independent over $\mathbb{R}$,
which implies that $B,B^2,\cdots,B^n$ are strongly independent over $\mathbb{Z}\setminus\{0\}$.\\
\rightline{$\Box$}

Let $f(n)$ denote det$(M_n(a))$ and $g(n)$ denote det$(N_n(a))$.
\begin{lemma}
\begin{equation*}
f(n)=\left\{
\begin{array}{cl}
1,&a=2\\
(-1)^{n-1}(2n+1),&a=-2\\
\frac{\alpha^{n+1}+\beta^{n+1}+\alpha^{n}+\beta^{n}}{2+a},& a \neq \pm 2
\end{array}\right.,~
g(n)=\left\{
\begin{array}{cl}
n+1,&a=2\\
(-1)^{n}(n+1),&a=-2\\
\frac{2(\alpha^{n}+\beta^{n})-a(\alpha^{n+1}+\beta^{n+1})}{4-a^2},& a \neq \pm 2
\end{array}\right.
\end{equation*}
where $\alpha,\beta$ are the roots of $\lambda^2-a\lambda+1=0$. i.e.
\begin{displaymath}
\alpha,\beta=\frac{a\pm\sqrt{a^2-4}}{2}.
\end{displaymath}
\end{lemma}
\noindent\textit{Proof}. Expanding the first columns we can get
\begin{equation}
f(n)=af(n-1)-f(n-2) ~~~\text{with}~~ f(1)=a-1~~\text{and}~~ f(2)=a^2-a-1.
\end{equation}
\begin{equation}
g(n)=ag(n-1)-g(n-2) ~~~\text{with}~~ g(1)=a~~\text{and}~~ g(2)=a^2-1.
\end{equation}
The characteristic equation of both is $\lambda^2-a\lambda+1=0$. Solve the difference equations we can get the results.\\
\rightline{$\Box$}

\begin{lemma}
The characteristic polynomial of $M_n(2)$ has $n$ different real roots:
\begin{equation*}
2+2\cos\left(\frac{2\pi}{2n+1}\right),2+2\cos\left(\frac{4\pi}{2n+1}\right),\cdots,2+2\cos\left(\frac{2n\pi}{2n+1}\right).
\end{equation*}
 Moreover
\begin{itemize}
  \item [(1)] all roots are irrational when $n\neq 1 \mod 3$,
  \item [(2)] 1 is the only rational root when $n=1 \mod 3$.
\end{itemize}
\end{lemma}
\noindent\textit{Proof}. If det$(M_n(a))=0$ for some $a$, then $a\neq \pm2$ and we have $\alpha^{n+1}+\beta^{n+1}+\alpha^n+\beta^n=0$. According to
$\alpha\beta=1$, we have $(\alpha+1)(\alpha^{2n+1}+1)=0$. Since $a\neq -2$, $\alpha\neq -1$. We have
\begin{equation*}
\alpha=e^{i\pi\left(\frac{2k}{2n+1}-1\right)},~k=1,2,\cdots,2n.
\end{equation*}
Therefore
\begin{equation*}
a=\alpha+\beta=2\text{Re}(\alpha)=-2\cos\left(\frac{2k\pi}{2n+1}\right),~k=1,2,\cdots,2n.
\end{equation*}
There are $n$ different values $-2\cos\left(\frac{2\pi}{2n+1}\right),-2\cos\left(\frac{4\pi}{2n+1}\right),\cdots,-2\cos\left(\frac{2n\pi}{2n+1}\right)$.
Note that the characteristic polynomial of $M_n(2)$ is $p(\lambda)=$det$(M_n(2-\lambda))$, consequently, which has $n$ different roots
\begin{equation*}
2+2\cos\left(\frac{2\pi}{2n+1}\right),2+2\cos\left(\frac{4\pi}{2n+1}\right),\cdots,2+2\cos\left(\frac{2n\pi}{2n+1}\right).
\end{equation*}
This proves the first part.\\

Note that det$(M_n(2))=1$ . We have
\begin{displaymath}
p(\lambda)=(-1)^n\lambda^n+\cdots+1.
\end{displaymath}
If there is a rational root of $p(\lambda)=0$ it must be $\pm 1$. By Proposition2, one has
\begin{eqnarray*}
p(1)=\frac{1}{3}\left(e^{i\frac{n+1}{3}\pi}+e^{-i\frac{n+1}{3}\pi}+e^{i\frac{n}{3}\pi}+e^{-i\frac{n}{3}\pi}\right)
= \frac{4}{3}\cos\frac{2n+1}{6}\pi\cos\frac{\pi}{6}
\end{eqnarray*}
and
\begin{equation*}
p(-1)=\frac{1}{5}\left[\left(\frac{3+\sqrt{5}}{2}\right)^{n+1}+\left(\frac{3-\sqrt{5}}{2}\right)^{n+1}+\left(\frac{3+\sqrt{5}}{2}\right)^{n}
+\left(\frac{3-\sqrt{5}}{2}\right)^{n}\right]>0.
\end{equation*}
This shows the second part.\\
\rightline{$\Box$}

\noindent\textit{Proof of Theorem \ref{strongly independent matrix}}. Since $2n+1$ is a prime number, we have $n\neq 1\mod 3$. Then all eigenvalues of $M_n(2)$ are different irrationals. Let the eigenvalues of $M_n(2)$ be $\lambda_1<\lambda_2<\cdots<\lambda_n$, and $\beta_i$ be the eigenvector corresponding to $\lambda_i$.
$\beta_i^{\perp}$ denotes the space of vectors that are orthogonal to $\beta_i$. \\
 Fix $i\in\{1,2,\cdots,n\}$. Let $B$ represent the $n\times n$ matrix
\begin{equation*}
\begin{bmatrix}
2-\lambda_i & 1& 0&\cdots&0&0\\
1&2-\lambda_i&1&\cdots&0&0\\
0&1&2-\lambda_i&\cdots&0&0\\
\vdots&\vdots&\vdots&&\vdots&\vdots\\
0&0&0&\cdots&2-\lambda_i&1\\
0&0&0&\cdots&1&1-\lambda_i
\end{bmatrix}
\end{equation*}
Then we have $\beta_i B=(0,0,\cdots,0)$. Denote the rows of $B$ by $\alpha_1,\alpha_2,\cdots,\alpha_{n}$. We get that
$\alpha_1,\alpha_2,\cdots,\alpha_n\in\beta_i^{\perp}$ and
\begin{displaymath}
\beta_i^{\perp}=\text{Span}\{\alpha_1,\alpha_2,\cdots,\alpha_{n-1}\},
\end{displaymath}
since $\alpha_1,\alpha_2,\cdots,\alpha_{n-1}$ are linearly independent.
If there exist $(x_1,x_2,\cdots,x_{n-1})\in\mathbb{R}^{n-1}$ and $(k_1,k_2,\cdots,k_n)\in\mathbb{Z}\setminus\{0\}$ such that
\begin{displaymath}
x_1\alpha_1+x_2\alpha_2+\cdots+x_{n-1}\alpha_{n-1}=(k_1,k_2,\cdots,k_n)
\end{displaymath}
That is
\begin{equation*}
\begin{bmatrix}
2-\lambda_i & 1& 0&\cdots&0\\
1&2-\lambda_i&1&\cdots&0\\
0&1&2-\lambda_i&\cdots&0\\
\vdots&\vdots&\vdots&&\vdots\\
0&0&0&\cdots&2-\lambda_i\\
0&0&0&\cdots&1
\end{bmatrix}
\begin{bmatrix}
x_1\\x_2\\ \vdots\\x_{n-1}
\end{bmatrix}
=\begin{bmatrix}
k_1\\k_2\\ \vdots\\k_{n-1}\\k_n
\end{bmatrix}
\end{equation*}
So $x_{n-1}=k_n$ and we have
\begin{equation*}
\begin{bmatrix}
x_1\\x_2\\ \vdots\\x_{n-2}\\ k_{n}
\end{bmatrix}=
\begin{bmatrix}
2-\lambda_i & 1& 0&\cdots&0\\
1&2-\lambda_i&1&\cdots&0\\
0&1&2-\lambda_i&\cdots&0\\
\vdots&\vdots&\vdots&&\vdots\\
0&0&0&\cdots&2-\lambda_i
\end{bmatrix}^{-1}
\begin{bmatrix}
k_1\\k_2\\ \vdots\\k_{n-1}
\end{bmatrix}
=N_{n-1}^{-1}(2-\lambda_i)\begin{bmatrix}
k_1\\k_2\\ \vdots\\k_{n-1}
\end{bmatrix}
\end{equation*}
Then there is a polynomial $q(x)\in \mathbb{Z}[x]$ of degree $n-1$ such that $q(\lambda_i)=0$ which is a contradiction. Since $2n+1$ is a prime number, $\varphi(2n+1)=2n$. By Lemma \ref{cycl1} and Lemma \ref{cycl2}, we have
\begin{displaymath}
\left[\mathbb{Q}\left(\cos\frac{2\pi}{2n+1},\cdots,\cos\frac{2n\pi}{2n+1}\right): \mathbb{Q}\right]=\frac{\varphi(2n+1)}{2}=n,
\end{displaymath}
which implies that it is impossible that for each eigenvalue there is a polynomial $q(x)\in \mathbb{Z}[x]$ of degree $n-1$ such that $q(\lambda_i)=0$.\\
\rightline{$\Box$}

\section{Characterizations of Ergodic, Weakly Mixing and Strongly Mixing via Fourier Coefficients}
\begin{theorem}\label{fourier coefficient}
The following are true.
\begin{itemize}
  \item [(1)] A measure $\mu$ on $\mathbb{T}^n$ is an ergodic $\times A$-invariant measure if and only if
  \begin{displaymath}
  \lim_{n\rightarrow \infty}\frac{1}{|F_n|}\sum_{j\in F_n} \hat{\mu}(\vec{k}A^j+\vec{l})=
  \hat{\mu}(\vec{k})\hat{\mu}(\vec{l})
  \end{displaymath}
  for every F{\o}lner sequence $\Sigma=\{F_n\}_{n=1}^\infty$ in $\mathbb{N}$ and $\vec{k},\vec{l}$ in $\mathbb{Z}^n$.
  \item [(2)]   A measure $\mu$ on $\mathbb{T}^n$ is a weakly mixing $\times A$-invariant measure if and only if
  \begin{displaymath}
  \lim_{n\rightarrow \infty}\frac{1}{|F_n|}\sum_{j\in F_n} |\hat{\mu}(\vec{k}A^j+\vec{l})-\hat{\mu}(\vec{k})\hat{\mu}(\vec{l})|^2=0
  \end{displaymath}
  for every F{\o}lner sequence $\Sigma=\{F_n\}_{n=1}^\infty$ in $\mathbb{N}$ and $\vec{k},\vec{l}$ in $\mathbb{Z}^n$.
  \item [(3)]  A measure $\mu$ on $\mathbb{T}^n$ is a strongly mixing $\times A$-invariant measure if and only if
  \begin{displaymath}
  \lim_{j\rightarrow \infty}\hat{\mu}(\vec{k}A^j+\vec{l})=\hat{\mu}(\vec{k})\hat{\mu}(\vec{l})
  \end{displaymath}
  for all $\vec{k},\vec{l}$ in $\mathbb{Z}^n$.
\end{itemize}
\end{theorem}

\noindent \textsl{Proof}.
\begin{itemize}
  \item [(1)] Suppose $\mu$ is an ergodic $\times A$-invariant measure on $\mathbb{T}^n$ and $T_A$ denotes the $\times A$ map.
  Using Lemma \ref{ergodicmean} for the measurable dynamical system $(\mathbb{T}^n,T_A,\mu)$ we have
  \begin{eqnarray}
  \lim_{n\rightarrow \infty} \frac{1}{|F_n|}\sum_{j\in F_n}\int_{\mathbb{T}^n} f(T_A^j x)g(x) d\mu(x)=\int_{\mathbb{T}^n} f d\mu \int_{\mathbb{T}^n} g d \mu\label{eq 5.1}
  \end{eqnarray}
  for every continuous functions $f,g$ on ${\mathbb{T}^n}$. Let $f=z^{\vec{k}}$ and $g=z^{\vec{l}}$ and we get
  \begin{eqnarray}
  \lim_{n\rightarrow \infty}\frac{1}{|F_n|}\sum_{j\in F_n} \hat{\mu}(\vec{k}A^j+\vec{l})=
  \hat{\mu}(\vec{k})\hat{\mu}(\vec{l})\label{eq 5.2}
  \end{eqnarray}
  which is the necessity.\\[1cm]
  Now we assume that (\ref{eq 5.1}) for every F{\o}lner sequence$\{F_n\}_{n=1}^\infty$ in $\mathbb{N}$ and all $\vec{k},\vec{l}\in\mathbb{Z}^n$. Let $\vec{l}=\vec{0}$, we get
  \begin{displaymath} \lim_{n\rightarrow \infty}\frac{1}{|F_n|}\sum_{j\in F_n} \hat{\mu}(\vec{k}A^j)=
  \hat{\mu}(\vec{k})\end{displaymath}
  for every $\vec{k}\in\mathbb{Z}^n$. Substituting $\vec{k}$ with $\vec{k}A$, we have
  \begin{eqnarray*}
  \hat{\mu}(\vec{k}A)
  &=& \lim_{n\rightarrow \infty}\frac{1}{|F_n|}\sum_{j\in F_n} \hat{\mu}(\vec{k}A^{j+1})\\
  &=& \lim_{n\rightarrow \infty}\frac{1}{|F_n+1|}\sum_{j\in F_n+1} \hat{\mu}(\vec{k}A^j)\\
  &=& \hat{\mu}(\vec{k}) ~~(\{F_n+1\}_{n=1}^\infty \text{is also a F{\o}lner sequence})
  \end{eqnarray*}
  for every $\vec{k}\in\mathbb{Z}^n$.  Hence $\mu$ is $\times A$-invariant.\\
  With (\ref{eq 5.2}) we have (\ref{eq 5.1}) is true for all  $f=z^{\vec{k}}$ and $g=z^{\vec{l}}$ . By linearity we have (\ref{eq 5.1}) is also true for all polynomials on $\mathbb{T}^n$. Since polynomials are dense in $L^2(\mathbb{T}^n,\mu)$, we have(\ref{eq 5.1})  true for all $f,g\in L^2(\mathbb{T}^n,\mu)$. Taking $f=g=1_E$ in (\ref{eq 5.1}) for some $T_A$-invariant Borel subset $E$ which satisfies $T^{-1}_A E=E$, we have $\mu(E)^2=\mu(E)$ . Hence $\mu$ is ergodic.

  \item [(2)] Suppose $\mu$ is a weakly mixing $\times A$-invariant measure on $\mathbb{T}^n$, which means $\mu\times\mu$ is an ergodic $T_A\times T_A$-invariant measure on $\mathbb{T}^{2n}$. Taking $f(z_1,z_2)=z_1^{\vec{k}}z_2^{\vec{k}}$ and $g(z_1,z_2)=z_1^{\vec{l}}z_2^{\vec{l}}$ in (\ref{eq ergodic mean}) of Lemma \ref{ergodicmean} with $X=\mathbb{T}^{2n}$ and $\nu=\mu\times\mu$, we get
      \begin{eqnarray}
      \lim_{n\rightarrow\infty}\sum_{j\in F_n}|\hat{\mu}(\vec{k}A^j+\vec{l})|^2=|\hat{\mu}(\vec{k})|^2|\hat{\mu}(l)|^2.\label{eq 5.3}
      \end{eqnarray}
      Since
      \begin{eqnarray*}
      &&|\hat{\mu}(\vec{k}A^j+\vec{l})-\hat{\mu}(\vec{k})\hat{\mu}(\vec{l})|^2\\
      &=& |\hat{\mu}(\vec{k}A^j+\vec{l})|^2+|\hat{\mu}(\vec{k})|^2|\hat{\mu}(l)|^2-\hat{\mu}
      (\vec{k}A^j+\vec{l})\hat{\mu}(-\vec{k})\hat{\mu}(-\vec{l})-
      \hat{\mu}
      (-\vec{k}A^j-\vec{l})\hat{\mu}(\vec{k})\hat{\mu}(\vec{l})
      \end{eqnarray*}
      We have
      \begin{eqnarray*}
      &&\lim_{n\rightarrow}\frac{1}{|F_n|}\sum_{j\in F_n}|\hat{\mu}(\vec{k}A^j+\vec{l})-\hat{\mu}(\vec{k})\hat{\mu}(\vec{l})|^2\\
      &=&\lim_{n\rightarrow}\frac{1}{|F_n|}\sum_{j\in F_n}|\hat{\mu}(\vec{k}A^j+\vec{l})|^2+|\hat{\mu}(\vec{k})|^2|\hat{\mu}(\vec{l})|^2-\hat{\mu}
      (\vec{k}A^j+\vec{l})\hat{\mu}(-\vec{k})\hat{\mu}(-\vec{l})-
      \hat{\mu}
      (-\vec{k}A^j-\vec{l})\hat{\mu}(\vec{k})\hat{\mu}(\vec{l})\\
      &=&|\hat{\mu}(\vec{k})|^2|\hat{\mu}(\vec{l})|^2+|\hat{\mu}(\vec{k})|^2|\hat{\mu}(\vec{l})|^2-|\hat{\mu}(\vec{k})|^2|\hat{\mu}(\vec{l})|^2-
      |\hat{\mu}(\vec{k})|^2|\hat{\mu}(\vec{l})|^2\\
      &=&0
      \end{eqnarray*}
      for all $\vec{k},\vec{l}\in\mathbb{Z}^n$ and every F{\o}lner sequence $\{F_n\}_{n=1}^\infty$ in $\mathbb{N}$. The second equality above is implied by (\ref{eq 5.3}) and the results of (1) noting that weakly mixing implies the ergodicity.\\
      Conversely, suppose
      \begin{eqnarray}
      \lim_{n\rightarrow}\frac{1}{|F_n|}\sum_{j\in F_n}|\hat{\mu}(\vec{k}A^j+\vec{l})-\hat{\mu}(\vec{k})\hat{\mu}(\vec{l})|^2=0\label{eq 5.4}
      \end{eqnarray}
      for all $\vec{k},\vec{l}\in\mathbb{Z}^n$ and every F{\o}lner sequence $\{F_n\}_{n=1}^\infty$ in $\mathbb{N}$.\\
      Firstly, (\ref{eq 5.4}) implies
       \begin{eqnarray*}
      \lim_{n\rightarrow}\frac{1}{|F_n|}\sum_{j\in F_n}\hat{\mu}(\vec{k}A^j+\vec{l})=\hat{\mu}(\vec{k})\hat{\mu}(\vec{l}).
      \end{eqnarray*}
      So by (1) $\mu$ is an ergodic $\times A$-invariant measure.\\
      In order to prove $\mu\times \mu$ is an erogodic $T_A\times T_A$-invariant measure on $\mathbb{T}^n$, it suffices to show that
      \begin{eqnarray*}
      \lim_{n\rightarrow \infty}\frac{1}{|F_n|}\sum_{j\in F_n} \int_{\mathbb{T}^n} f((T_A\times T_A)^j(z_1,z_2))g(z_1,z_2)d\mu\times\mu(z_1,z_2)=\int_{\mathbb{T}^n}fd\mu
      \int_{\mathbb{T}^n}g d\mu
      \end{eqnarray*}
      for all continuous functions $f$ and $g$ on $\mathbb{T}^n$, which is equivalent to show
      \begin{eqnarray*}
       \lim_{n\rightarrow \infty}\frac{1}{|F_n|}\sum_{j\in F_n} \hat{\mu}(\vec{k_1}A^j+\vec{l_1})\hat{\mu}(\vec{k_2}A^j+\vec{l_2})=
       \hat{\mu}(\vec{k_1})\hat{\mu}(\vec{k_2})\hat{\mu}(\vec{l_1})\hat{\mu}(\vec{l_2})
       \end{eqnarray*}
       for all $\vec{k_1},\vec{k_2},\vec{l_1},\vec{l_2}\in \mathbb{Z}^n$ by taking $f=z_1^{\vec{k_1}}z_2^{\vec{k_2}}$ and $g=z_1^{\vec{l_1}}z_2^{\vec{l_2}}$ since whose linear span is dense in $C(\mathbb{T}^{2n})$.\\
       Note that
       \begin{eqnarray*}
       &&\left|\hat{\mu}(\vec{k_1}A^j+\vec{l_1})\hat{\mu}(\vec{k_2}A^j+\vec{l_2})-
       \hat{\mu}(\vec{k_1})\hat{\mu}(\vec{k_2})\hat{\mu}(\vec{l_1})\hat{\mu}(\vec{l_2})\right|\\
       &\leq&  \left|\hat{\mu}(\vec{k_1}A^j+\vec{l_1})\left[\hat{\mu}(\vec{k_2}A^j+\vec{l_2})-
       \hat{\mu}(\vec{k_2})(\vec{l_2})\right]\right|
       + \left|\left[\hat{\mu}(\vec{k_1}A^j+\vec{l_1})-\hat{\mu}(\vec{k_1})\hat{\mu}(\vec{l_1})\right]
       \hat{\mu}(\vec{k_2})(\vec{l_2})\right|\\
       &\leq& \left|\hat{\mu}(\vec{k_2}A^j+\vec{l_2})-\hat{\mu}(\vec{k_2})(\vec{l_2})\right|
        +\left|\hat{\mu}(\vec{k_1}A^j+\vec{l_1})-\hat{\mu}(\vec{k_1})\hat{\mu}(\vec{l_1})\right|
       \end{eqnarray*}
       for all $\vec{k_1},\vec{k_2},\vec{l_1},\vec{l_2}\in \mathbb{Z}^n$.\\
       Hence we have
       \begin{eqnarray*}
       &&\lim_{n\rightarrow \infty}\frac{1}{|F_n|}\sum_{j\in F_n} \left|\hat{\mu}(\vec{k_1}A^j+\vec{l_1})\hat{\mu}(\vec{k_2}A^j+\vec{l_2})-
       \hat{\mu}(\vec{k_1})\hat{\mu}(\vec{k_2})\hat{\mu}(\vec{l_1})\hat{\mu}(\vec{l_2})\right|^2\\
       &\leq& \lim_{n\rightarrow \infty}\frac{1}{|F_n|}\sum_{j\in F_n}\left[\left|\hat{\mu}(\vec{k_2}A^j+\vec{l_2})-\hat{\mu}(\vec{k_2})(\vec{l_2})\right|
        +\left|\hat{\mu}(\vec{k_1}A^j+\vec{l_1})-\hat{\mu}(\vec{k_1})\hat{\mu}(\vec{l_1})\right|\right]^2\\
        &\leq& 2\lim_{n\rightarrow \infty}\frac{1}{|F_n|}\sum_{j\in F_n}\left[\left|\hat{\mu}(\vec{k_2}A^j+\vec{l_2})-\hat{\mu}(\vec{k_2})(\vec{l_2})\right|^2
        +\left|\hat{\mu}(\vec{k_1}A^j+\vec{l_1})-\hat{\mu}(\vec{k_1})\hat{\mu}(\vec{l_1})\right|^2\right]\\
        &=&0~~(\text{by (\ref{eq 5.4})}).
        \end{eqnarray*}
        The last inequality above is implied by the inequality $(a+b)^2\leq 2(a^2+b^2)$ for any nonnegative $a,b$. Using Jensen's inequality
        $\left(\frac{|x_1|+\cdots+|x_n|}{n}\right)^2\leq \frac{|x_1|^2+\cdots+|x_n|^2}{n}$, we get
        \begin{displaymath}
        \lim_{n\rightarrow \infty}\frac{1}{|F_n|}\sum_{j\in F_n} \left|\hat{\mu}(\vec{k_1}A^j+\vec{l_1})\hat{\mu}(\vec{k_2}A^j+\vec{l_2})-
       \hat{\mu}(\vec{k_1})\hat{\mu}(\vec{k_2})\hat{\mu}(\vec{l_1})\hat{\mu}(\vec{l_2})\right|=0.
       \end{displaymath}
       This proves the sufficiency.
  \item [(3)] Suppose $\mu$ is strongly mixing, which means that $\lim_{j\rightarrow\infty}\mu(T_A^{-j}E\cap F)=\mu(E)\mu(F)$ for all Borel subsets $E,F$
  in $\mathbb{T}^n$. Let $1_E$ stand for the characteristic function of $E$, then we have
  \begin{eqnarray*}
  \lim_{j\rightarrow \infty} \int_{\mathbb{T}^n} 1_E(T_A^jx)1_F(x)d\mu(x)=\int_{\mathbb{T}^n}1_Ed\mu \int_{\mathbb{T}^n}1_Fd\mu
  \end{eqnarray*}
  for all Borel subsets $E,F$ in $\mathbb{T}^n$.\\
  Since the linear combinations of characteristic functions are dense in $L^2(\mathbb{T}^n,\mu)$, we have
  \begin{eqnarray*}
  \lim_{j\rightarrow \infty} \int_{\mathbb{T}^n} f(T_A^jx)g(x)d\mu(x)=\int_{\mathbb{T}^n}fd\mu \int_{\mathbb{T}^n}gd\mu
  \end{eqnarray*}
  for all $f,g\in C(\mathbb{T}^n)$. Particularly taking $f=z^{\vec{k}}$ and $g=z^{\vec{l}}$, we have
  \begin{eqnarray}
  \lim_{j\rightarrow\infty}\hat{\mu}(\vec{k}A^j+\vec{l})=\hat{\mu}(\vec{k})\hat{\mu}(\vec{l})\label{eq 5.5}
  \end{eqnarray}
  for all $\vec{k},\vec{l}$ in $\mathbb{Z}^n$. This proves the necessity.\\
  On the other hand, suppose a measure $\mu$ on $\mathbb{T}^n$ satisfies (\ref{eq 5.5}). Let $\vec{l}=\vec{0}$ and replace $\vec{k}$ by $\vec{k}A$. Then we get
  \begin{displaymath}
  \hat{\mu}(\vec{k}A)=\lim_{j\rightarrow \infty}\hat{\mu}(\vec{k}A^{j+1})=\lim_{j\rightarrow \infty}\hat{\mu}(\vec{k}A^{j})=\hat{\mu}(\vec{k})
  \end{displaymath}
  for all $\vec{k}\in\mathbb{Z}^n$. Hence $\mu$ is $\times A$-invariant and we have
  \begin{eqnarray*}
  \lim_{j\rightarrow \infty} \int_{\mathbb{T}^n} f(T_A^jx)g(x)d\mu(x)=\int_{\mathbb{T}^n}fd\mu \int_{\mathbb{T}^n}gd\mu
  \end{eqnarray*}
  when $f=z^{\vec{k}}$ and $g=z^{\vec{l}}$ for all $\vec{k},\vec{l}$ in $\mathbb{Z}^n$. Since the linear combinations of $z^{\vec{k}}$ and $z^{\vec{l}}$ are polynomials on $\mathbb{T}^n$ which are dense in $L^2(\mathbb{T}^n,\mu)$, the above is also true for all  $f,g\in L^2(\mathbb{T}^n,\mu)$. Particularly it holds for $f=1_E$ and $g=1_F$ for any Borel subsets $E,F$ in $\mathbb{T}^n$ that is
  \begin{displaymath}
  \lim_{j\rightarrow\infty}\mu(T_A^{-j}E\cap F)=\mu(E)\mu(F).
  \end{displaymath}
  Hence we complete the proof.
\end{itemize}
\rightline{$\Box$}
\section{Measure Rigidity on $\mathbb{T}^n$}
\textit{Proof of Theorem \ref{rigidity on torus}}.
\begin{itemize}
  \item [(1)] Suppose $\mu$ is an ergodic $\times A$-invariant measure on $\mathbb{T}^n$ and there exist a sequence of $n$ matrices $B_1, B_2, \cdots, B_n\in GL(n,\mathbb{Z})$ which are strongly independent over $\mathbb{Z}\setminus\{0\}$ and a F{\o}lner sequence $\Sigma=\{F_n\}_{n=1}^\infty$ such that $\mu$ is $\times(A^j+B_i)$-invariant for all $i\in\{1,2,\cdots,n\}$ and $j$ in some $E\subseteq \mathbb{N}$ with $D_{\Sigma}(E)=1$.\\[0.2cm]
      If $\mu$ is not Lebesgue measure, then there exists a nonzero $\vec{k}\in \mathbb{Z}^n$ such that $\hat{\mu}(\vec{k})\neq 0$.\\[0.2cm]
      Since $\mu$ is an ergodic $\times A$-invariant measure, by Theorem \ref{fourier coefficient} (1), we have
      \begin{eqnarray*}
       \lim_{n\rightarrow \infty}\frac{1}{|F_n|}\sum_{j\in F_n} \hat{\mu}(\vec{k}A^j+\vec{k}B_i)=
       \hat{\mu}(\vec{k})\hat{\mu}(\vec{k}B_i).
      \end{eqnarray*}
      for every $i\in\{1,2,\cdots,n\}$. Note that
      \begin{eqnarray*}
      \frac{1}{|F_n|}\sum_{j\in F_n} \hat{\mu}(\vec{k}A^j+\vec{k}B_i)
      &=& \frac{1}{|F_n|}\left(\sum_{j\in F_n\cap E}+\sum_{j\in F_n\setminus E}\right) \hat{\mu}(\vec{k}A^j+\vec{k}B_i)\\
      &=& \frac{|F_n\cap E|}{|F_n|}\hat{\mu}(\vec{k})+\sum_{j\in F_n\setminus E} \hat{\mu}(\vec{k}A^j+\vec{k}B_i)\\
      &\rightarrow&\hat{\mu}(\vec{k})~~~~\text{as} ~~ n\rightarrow \infty.
      \end{eqnarray*}
      Hence $\hat{\mu}(\vec{k})=\hat{\mu}(\vec{k})\hat{\mu}(\vec{k}B_i)$ which implies $\hat{\mu}(\vec{k}B_i)=1$ for every $i\in\{1,2,\cdots,n\}$. According to Lemma \ref{atomic}, we get a contradiction.
  \item [(2)] Suppose $\mu$ is an weakly mixing $\times A$-invariant measure on $\mathbb{T}^n$ and there exist a sequence of $n$ matrices $B_1, B_2, \cdots, B_n\in GL(n,\mathbb{Z})$ which are strongly independent over $\mathbb{Z}\setminus\{0\}$ and a F{\o}lner sequence $\Sigma=\{F_n\}_{n=1}^\infty$ such that $\mu$ is $\times(A^j+B_i)$-invariant for all $i\in\{1,2,\cdots,n\}$ and $j$ in some $E\subseteq \mathbb{N}$ with $\overline{D}_{\Sigma}(E)>0$.\\[0.2cm]
      If $\mu$ is not Lebesgue measure, then there exists a nonzero $\vec{k}\in \mathbb{Z}^n$ such that $\hat{\mu}(\vec{k})\neq 0$.\\[0.2cm]
      Since $\mu$ is a weakly mixing $\times A$-invariant measure, by Theorem \ref{fourier coefficient} (2), we have
      \begin{eqnarray*}
       \lim_{n\rightarrow \infty}\frac{1}{|F_n|}\sum_{j\in F_n} |\hat{\mu}(\vec{k}A^j+\vec{k}B_i)-\hat{\mu}(\vec{k})\hat{\mu}(\vec{k}B_i)|^2=0.
      \end{eqnarray*}
      Therefore,
      \begin{eqnarray*}
      0&=&  \limsup_{n\rightarrow \infty}\frac{1}{|F_n|}\sum_{j\in F_n} |\hat{\mu}(\vec{k}A^j+\vec{k}B_i)-\hat{\mu}(\vec{k})\hat{\mu}(\vec{k}B_i)|^2\\
      &\geq& \limsup_{n\rightarrow \infty}\frac{1}{|F_n|}\sum_{j\in F_n\cap E} |\hat{\mu}(\vec{k}A^j+\vec{k}B_i)-\hat{\mu}(\vec{k})\hat{\mu}(\vec{k}B_i)|^2\\
      &=& \limsup_{n\rightarrow \infty}\frac{1}{|F_n|}\sum_{j\in F_n\cap E} \left[ |\hat{\mu}(\vec{k}A^j+\vec{k}B_i)|^2+|\hat{\mu}(\vec{k})|^2|\hat{\mu}(\vec{k}B_i)|^2\right.\\
      &&\qquad\qquad\quad \left. -\hat{\mu}(\vec{k}A^j+\vec{k}B_i)\hat{\mu}(-\vec{k})\hat{\mu}(-\vec{k}B_i)-\hat{\mu}(-\vec{k}A^j-\vec{k}B_i)\hat{\mu}(\vec{k})\hat{\mu}(\vec{k}B_i)\right]\\
      &=& \limsup_{n\rightarrow \infty}\frac{1}{|F_n|}\sum_{j\in F_n\cap E}\left[ |\hat{\mu}(\vec{k})|^2+|\hat{\mu}(\vec{k})|^2|\hat{\mu}(\vec{k}B_i)|^2
     -|\hat{\mu}(\vec{k})|^2\hat{\mu}(-\vec{k}B_i)-|\hat{\mu}(\vec{k})|^2\hat{\mu}(\vec{k}B_i)\right]\\
     &=& \limsup_{n\rightarrow \infty}\frac{1}{|F_n|}\sum_{j\in F_n\cap E} |\hat{\mu}(\vec{k})-\hat{\mu}(\vec{k})\hat{\mu}(\vec{k}B_i)|^2\\
     &=& |\hat{\mu}(\vec{k})-\hat{\mu}(\vec{k})\hat{\mu}(\vec{k}B_i)|^2 \overline{D}_{\Sigma}(E).
      \end{eqnarray*}
      Hence $\hat{\mu}(\vec{k})-\hat{\mu}(\vec{k})\hat{\mu}(\vec{k}B_i)=0$ which implies that $\hat{\mu}(\vec{k}B_i)=1$ for every $i\in\{1,2,\cdots,n\}$. We get a contradiction again.
  \item [(3)]  Suppose $\mu$ is a strongly mixing $\times A$-invariant measure on $\mathbb{T}^n$ and there exist a sequence of $n$ matrices $B_1, B_2, \cdots, B_n\in GL(n,\mathbb{Z})$ which are strongly independent over $\mathbb{Z}\setminus\{0\}$ and an infinite set $E\subseteq\mathbb{N}$ such that $\mu$ is $\times(A^j+B_i)$-invariant for all $i\in\{1,2,\cdots,n\}$ and $j$ in $E$.\\[0.2cm]
      If $\mu$ is not Lebesgue measure, then there exists a nonzero $\vec{k}\in \mathbb{Z}^n$ such that $\hat{\mu}(\vec{k})\neq 0$.\\[0.2cm]
      Since $\mu$ is a strongly mixing $\times A$-invariant measure, by Theorem \ref{fourier coefficient} (3), we have
      \begin{eqnarray*}
      \lim_{j\rightarrow \infty\atop j\in E} \hat{\mu}(\vec{k}A^j+\vec{k}B_i)=\hat{\mu}(\vec{k})\hat{\mu}(\vec{k}B_i)
      \end{eqnarray*}
      for all $i\in\{1,2,\cdots,n\}$ and $j$ in $E$. Owing to $\mu$ being $\times(A^j+B_i)$-invariant for all $j$ in $E$, one has $\hat{\mu}(\vec{k}A^j+\vec{k}B_i)=\hat{\mu}(\vec{k})$. Consequently, $\hat{\mu}{\vec{k}}=\hat{\mu}(\vec{k})\hat{\mu}(\vec{k}B_i)$, which implies $\hat{\mu}(\vec{k}B_i)=1$ for all $i\in\{1,2,\cdots,n\}$. This again leads to a contradiction.
\end{itemize}
\indent Suppose $\mu$ is a measure on $\mathbb{T}^n$ satisfying (2) or (3) of Theorem \ref{rigidity on torus}. If $\mu$ is not a Lebesgue measure, then $\mu$ is atomic. According to Lemma \ref{dirac}, we claim that $\mu$ is a Dirac measure on $\mathbb{T}^n$.\\
\rightline{$\Box$}

\vspace{3mm}

{\noindent\bf Acknowledgement.} We would like to thank professor Shengkui Ye and professor Yi Gu for their helpful discussions. The second author is supported by NSFC (No. 11771318, No. 11790274).

\addcontentsline{toc}{section}{Reference}
\bibliography{bib}

\begin{thebibliography}{99}
\bibitem{F}H. Furstenberg, Disjointness in ergodic theory, minimal sets, and a problem in Diophantine approximation,\emph{ Math. Systems Theory}, 1 (1967), 1-49.
\bibitem{H}H. Huang, Fourier Coefficients of $\times p$-invariant Measures, \emph{Journal of Modern Dynamics}, 11 (2017), 551-562.
\bibitem{KK}B. Kalinin and A. Katok, Invariant measures for actions of higher rank abelian groups, Proceedings of Symposia in Pure Mathematics, 2009, 593-637.
\bibitem{KS1} A. Katok and R. Spatzier, Invariant Measures for Higher Rank Hyperbolic Abelian Actions, \emph{Ergod. Th. and Dynam. Syst.}, 16 (1996), 751-778.
\bibitem{KS2} A. Katok and R. Spatzier, Corrections to ``Invariant measures for higher rank hyperbolic abelian actions", \emph{Ergod. Th. and Dynam. Syst.}, 18 (1998), 503-507.
\bibitem{L} R. Lyons, On measures simultaneously 2- and 3-invariant, \emph{Israel J. Math.}, 61 (1988), 219-224.
\bibitem{R}D. Rudolph, $\times 2$ and $\times 3$ invariant measures and entropy, \emph{Ergod. Th. and Dynam. Syst.}, 10(1990), 395-406.
\bibitem{W} S. Weintraub., \emph{Galois Theorey}, Springer; 2nd ed. 2009 edition (November 21, 2008).

\end{thebibliography}

\vspace{5mm}
\noindent Huichi Huang, College of Mathematics and Statistics, Chongqing University, Chongqing, 401331, P.R. China (huanghuichi@cqu.edu.cn)\\

\noindent Enhui Shi, School of mathematical and sciences, Soochow University, Suzhou, 215006, P.R. China (ehshi@suda.edu.cn)\\

\noindent Hui Xu, School of mathematical and sciences, Soochow University, Suzhou, 215006, P.R. China (mathegoer@163.com)

\end{document}